\documentclass[a4paper,11pt]{scrartcl}

\usepackage[T1]{fontenc}
\usepackage[utf8]{inputenc}
\usepackage[english]{babel}
\usepackage{lmodern}
\usepackage{microtype}

\usepackage{mleftright}
\usepackage[inline,shortlabels]{enumitem}
\setlist{font=\normalfont}
\usepackage{graphicx}
\usepackage{amsmath}
\usepackage{amssymb}
\usepackage{amsthm}
\usepackage{tikz}
\usepackage[noadjust]{cite}
\usepackage[hidelinks]{hyperref}
\usepackage[capitalize,noabbrev]{cleveref}

\theoremstyle{definition}
\newtheorem{definition}{Definition}[section]

\newtheorem{example}[definition]{Example}
\newtheorem{remark}[definition]{Remark}

\theoremstyle{plain}

\newtheorem{proposition}[definition]{Proposition}

\newcommand{\R}{\mathbb{R}}
\newcommand{\N}{\mathbb{N}}
\newcommand{\wto}{\rightharpoonup}

\title{
    Brezis pseudomonotonicity is strictly weaker than Ky--Fan hemicontinuity%
    \thanks{This research was supported by the German Research Foundation (DFG) within the priority program ``Non-smooth and Complementarity-based Distributed Parameter Systems: Simulation and Hierarchical Optimization'' (SPP 1962) under grant number KA 1296/24-1.}
}

\date{October 9, 2018}

\author{
    Daniel Steck%
    \thanks{University of W\"urzburg, Institute of Mathematics, Emil-Fischer-Str.\ 30, 97074 Würzburg, Germany; email: \href{mailto:mail@danielsteck.net}{\nolinkurl{mail@danielsteck.net}}. ORCID: \href{https://orcid.org/0000-0001-5335-5428}{\nolinkurl{0000-0001-5335-5428}}.}
}

\begin{document}

\maketitle

{
\small{\bfseries\sffamily\abstractname.}
    In 1968, H.~Brezis introduced a notion of operator pseudomonotonicity which provides a unified approach to monotone and nonmonotone variational inequalities (VIs). A closely related notion is that of Ky--Fan hemicontinuity, a continuity property which arises if the famous Ky--Fan minimax inequality is applied to the VI framework. It is clear from the corresponding definitions that Ky--Fan hemicontinuity implies Brezis pseudomonotonicity, but quite surprisingly, a recent publication by Sadeqi and Paydar (J.\ Optim.\ Theory Appl., 165(2):344--358, 2015) claims the equivalence of the two properties. The purpose of the present note is to show that this equivalence is false; this is achieved by providing a concrete example of a nonlinear operator which is Brezis pseudomonotone but \emph{not} Ky--Fan hemicontinuous.
\par\addvspace{\baselineskip}
}

{
\par\addvspace{\baselineskip}
}

{
\par\addvspace{\baselineskip}
}

\section{Introduction}

Variational inequalities (VIs) are a prominent tool in applied mathematics. They have found numerous applications, including constrained optimization problems, Nash equilibrium problems, and several types of contact problems in mechanics. More details can be found, for instance, in \cite{Facchinei2003,Kinderlehrer2000,Baiocchi1984,Glowinski1981} and the references therein.

The study of variational inequalities can be divided into multiple facets: most commonly, one is interested in sufficient conditions for the existence of solutions, the design of suitable algorithms for their computation, or other properties of the solution set such as closedness or convexity. Throughout the last decades, various concepts have been developed in order to ascertain these properties, including the monotonicity of the variational operator (and multiple relaxed versions thereof) as well as several types of continuity (hemicontinuity, continuity on finite-dimensional subspaces, etc.). In addition, one often requires suitable properties of the feasible set such as closedness or (weak) compactness.

One of the most general properties which can be used to tackle variational inequalities is that of \emph{(Brezis) pseudomonotonicity}, a property which was introduced in \cite{Brezis1968}. (This should not be confused with pseudomonotonicity in the sense of Karamardian.) The attractive property of Brezis pseudomonotonicity is that it provides a unified approach to monotone and nonmonotone problems---indeed, it is best viewed as a hybrid combining elements of both monotonicity and continuity, see \cref{Dfn:BrezisPseudomonotonicity}.

Since its conception in 1968, the notion of pseudomonotonicity has occurred prominently in the works of Browder \cite{Browder1972}; Brezis, Nirenberg, and Stampacchia \cite{Brezis1972}; Zeidler \cite{Zeidler1990}; and Barbu and Precupanu \cite{Barbu2012}. The standard application of pseudomonotonicity was the construction of existence results for VIs, a topic which occurs in all these references and was also revisited in \cite{Kien2009}. In addition, pseudomonotonicity has turned out to be quite useful when analyzing convergence of iterative algorithms for constrained minimization and variational or quasi-variational inequalities, see \cite{Kanzow2018,Steck2018} for more details.

A different but related approach to the existence of solutions is given by the classical minimax inequality of Ky Fan \cite{Fan1972}. An application of this result to the VI framework gives rise to a continuity property which is sometimes called \emph{Ky--Fan hemicontinuity}. This property implies Brezis pseudomonotonicity (a fact which follows directly from the corresponding definitions, see below), but the latter appears to be more refined and convenient when dealing with infinite-dimensional VIs. However, quite surprisingly, a 2015 publication by Sadeqi and Paydar \cite{Sadeqi2015} claims the equivalence of the two properties. The purpose of the present paper is to discuss this equivalence and provide a counterexample which shows that the two properties are in fact distinct. In addition, we also outline an error in the reference which may have led to the false result.

This paper is organized as follows. In \cref{Sec:Prelims}, we give a brief summary of the properties in question, their consequences, and relations to other standard properties for VIs. \Cref{Sec:Example} contains the main counterexample.

\section{Hemicontinuity, pseudomonotonicity, and their role in the study of variational inequalities}\label{Sec:Prelims}

Throughout this paper, $X$ is a real Banach space with norm $\|\cdot\|_X$ and continuous dual $X^*$. The duality pairing between $X^*$ and $X$ is denoted by $\langle\cdot,\cdot\rangle$. We write $\to$, $\wto$, and $\wto^*$ for strong, weak, and weak-$^*$ convergence.

Let $A\subseteq X$ be a nonempty convex set and $F:X\to X^*$ a given mapping. The typical variational inequality involving $F$ and $A$ takes on the following form:
\begin{equation}\label{Eq:VI}
    \text{Find }\bar{x}\in A\text{ such that}\quad
    \langle F(\bar{x}),y-\bar{x} \rangle \ge 0
    \quad\text{for all }y\in A.
\end{equation}
A simple but useful observation is that the VI \eqref{Eq:VI} can be rewritten as an \emph{equilibrium problem} in the following sense: define $\Psi:A^2\to\R$, $\Psi(x,y):=\langle F(x),x-y \rangle$. Then the VI is obviously equivalent to the existence of $\bar{x}\in A$ such that $\Psi(\bar{x},y)\le 0$ for all $y\in A$. One of the standard existence results for such problems is the minimax inequality of Ky Fan \cite{Fan1972}, which requires the weak sequential lower semiconinuity of $\Psi$ with respect to $x$. This gives rise to the following definition.

\begin{definition}[Ky--Fan hemicontinuity]\label{Dfn:KyFanHemicontinuity}
    We say that $F:X\to X^*$ is \emph{Ky--Fan hemicontinuous} if, for every $y\in X$, the function $x\mapsto \langle F(x),x-y \rangle$ is weakly sequentially lower semicontinuous.
\end{definition}

The above is one of the two main properties which we will discuss in this paper. The second one was introduced by Brezis \cite{Brezis1968} and is given as follows.

\begin{definition}[Brezis pseudomonotonicity]\label{Dfn:BrezisPseudomonotonicity}
    We say that an operator $F:X\to X^*$ is \emph{(Brezis) pseudomonotone} if, whenever
\begin{equation}\label{Eq:DfnBrezisPseudomonotonicity1}
    \{x^k\}\subseteq X, \quad x^k\wto x, \quad\text{and}\quad
    \limsup_{k\to\infty}{}\langle F(x^k),x^k-x \rangle \le 0,
\end{equation}
    then
\begin{equation}\label{Eq:DfnBrezisPseudomonotonicity2}
    \langle F(x),x-y \rangle \le \liminf_{k\to\infty}{} \langle F(x^k),x^k-y \rangle
    \quad\text{for all }y\in X.
\end{equation}
\end{definition}

It is clear from the above definitions that Brezis pseudomonotonicity is weaker than Ky--Fan hemicontinuity: the latter requires that \eqref{Eq:DfnBrezisPseudomonotonicity2} holds for all weakly convergent sequences $\{x^k\}\subseteq X$ with limit $x\in X$, whereas Brezis pseudomonotonicity only asserts this estimate for sequences which additionally satisfy the $\limsup$-condition in \eqref{Eq:DfnBrezisPseudomonotonicity1}.

The set of pseudomonotone operators is large and encompasses many practically relevant examples. Various sufficient conditions for pseudomonotonicity can be found in \cite{Zeidler1990,Steck2018}; in particular, $F$ is pseudomonotone provided it is either
\begin{enumerate*}[(i)]
    \item monotone and continuous,
    \item completely continuous, or
    \item the sum of two operators which are themselves pseudomonotone.
\end{enumerate*}
Using results from differential calculus in Banach spaces, one can also give sufficient conditions for pseudomonotonicity in the special case where $F$ is the Fr\'echet-derivative of a real-valued functional, see \cite{Steck2018}.

The following is the basic existence result for VIs with pseudomonotone operators. Note that we call $F$ \emph{bounded} if it maps bounded sets in $X$ to bounded sets in $X^*$.

\begin{proposition}[Pseudomonotone VIs, {\cite[Corollary~4.2]{Kanzow2018}}]\label{Prop:ExistencePseudo}
    Let $A\subseteq X$ be a nonempty, convex, weakly compact set, and $F:X\to X^*$ a bounded pseudomonotone operator. Then the variational inequality \eqref{Eq:VI} admits a solution $\hat{x}\in A$.
\end{proposition}


Apart from the existence of solutions, the notion of pseudomonotonicity has some additional interesting consequences. Two particular ones are mentioned below:
\begin{itemize}
\item If $F:X\to X^*$ is pseudomonotone, then the solution set of the VI \eqref{Eq:VI} is always weakly sequentially closed. More generally, if $\{x^k\}$ is a sequence of suitable ``approximate'' solutions of the VI, then every weak limit point of $\{x^k\}$ belongs to its solution set, see \cite{Kanzow2018}.
\item If $X$ is \emph{finite-dimensional} (without loss of generality, a Hilbert space), then an operator $F:X\to X$ is bounded and pseudomonotone if and only if it is continuous. Thus, in this case, the study of pseudomonotone VIs subsumes the well-known theory of finite-dimensional VIs; see, for instance, \cite{Facchinei2003}.
\end{itemize}
We close this section with a general comment on the role of sequences in \cref{Dfn:KyFanHemicontinuity,Dfn:BrezisPseudomonotonicity}. This remark is particularly important when using the analogues of these properties in topological vector spaces.

\begin{remark}[Sequences versus nets]\label{Rem:Nets}
    Some authors define the aforementioned concepts on a general Hausdorff topological vector space (instead of a Banach space endowed with its weak topology). In that case, the properties ought to be formulated in terms of nets or filters instead of ordinary sequences; see, for instance, \cite{Brezis1972}.
\end{remark}

\section{A nonlinear operator which is Brezis pseudomonotone but not Ky--Fan hemicontinuous}\label{Sec:Example}

As we shall see below, Brezis pseudomonotonicity and Ky--Fan hemicontinuity are not equivalent. The particular example shown here involves the well-known function spaces $L^p(\Omega)$, $W^{k,p}(\Omega)$, and $W_0^{k,p}(\Omega)$, with $\Omega$ a bounded finite-dimensional domain, $k\in\N$, and $p\in [1,+\infty]$; see, for instance, \cite{Adams2003}.

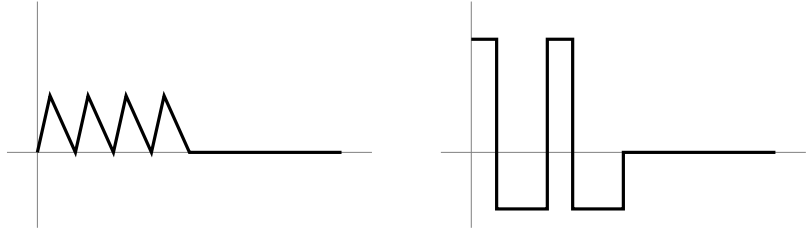
\begin{figure}\centering
\begin{tikzpicture}[xscale=4,yscale=3]
    \draw[very thin,color=gray] (-0.1,0) -- (1.1,0); 
    \draw[very thin,color=gray] (0,-1/3) -- (0,2/3); 
    
    \draw[very thick] (0,0) -- (1/24,1/4) -- (1/8,0) -- (4/24,1/4) -- (2/8,0) -- (7/24,1/4) -- (3/8,0) -- (10/24,1/4) -- (4/8,0) -- (1,0);
\end{tikzpicture}
\qquad
\begin{tikzpicture}[xscale=4,yscale=0.25]
    \draw[very thin,color=gray] (-0.1,0) -- (1.1,0); 
    \draw[very thin,color=gray] (0,-4) -- (0,8); 
    
    \draw[very thick] (0,6) -- (1/12,6) -- (1/12,-3) -- (1/4,-3) -- (1/4,6) -- (4/12,6) -- (4/12,-3) -- (2/4,-3) -- (2/4,0) -- (4/4,0);
\end{tikzpicture}
\caption{The sequence from \cref{Ex:PseudoKyFan} for $k=4$: $u_k$ (left) and $\nabla u_k$ (right).}
\label{Fig:ExPseudoKyFan}
\end{figure}

\begin{example}\label{Ex:PseudoKyFan}
    Let $X:=W_0^{1,3} (0,1)$ be a well-known Sobolev space on the interval $(0,1)$, and let $F:X\to X^*$ be the (negative) $p$-Laplacian defined by
\begin{equation*}
    \langle F(u),v \rangle:= \int_0^1 |\nabla u(t)|
    \nabla u(t) \nabla v(t) \,\textnormal{d}t.
\end{equation*}
    More details on this operator can be found in \cite[Section~3.1.3]{Ruzicka2006}. In particular, it is known that $F$ is monotone and continuous, hence pseudomonotone. Now, for each $k\in\N$, let $u_k:[0,1]\to\R$ be the piecewise linear function with value $1/k$ at $t=(3 i+1)/(6 k)$ for $i=0,\ldots,k-1$, and value zero at $t=i/(2 k)$ for $i=0,\ldots,k$, and on $[1/2,1]$. Clearly, $u_k\to 0$ in $L^3(0,1)$. Moreover, the weak derivative of $u_k$ is (almost everywhere) given by
\begin{equation*}
    \nabla u_k(t)=
\begin{cases}
    6, & \text{if }t\in \bigl( \frac{i}{2 k},\frac{3 i+1}{6 k} \bigr)
    \text{ with }i=0,\ldots,k-1, \\
    -3, & \text{if }t\in \bigl( \frac{3 i-2}{6 k},\frac{i}{2 k} \bigr)
    \text{ with }i=1,\ldots,k, \\
    0, & \text{if }t\in \bigl( \frac{1}{2},1 \bigr),
\end{cases}
\end{equation*}
    see \cref{Fig:ExPseudoKyFan}. It follows from standard arguments (e.g., \cite[Exercise~8.7,~p.\,254]{Alt2016}) that $\nabla u_k\wto 0$ in $L^3(0,1)$. Hence, $u_k\wto u:=0$ in $X=W_0^{1,3}(0,1)$. Now, let $v(t):=\alpha \min\{t,1-t\}$ with $\alpha\gg 0$, and observe that $\|\nabla u_k \|_{L^3(0,1)}^3=45$ for all $k$. Finally, we have $\langle F(u),u-v \rangle=0$, but an elementary calculation shows that
\begin{equation*}
    \langle F(u_k),u_k-v \rangle = \|\nabla u_k \|_{L^3(0,1)}^3-\int_0^1 |\nabla u_k(t)| \nabla u_k(t)
    \nabla v(t) \,\textnormal{d}t=45-3\alpha.
\end{equation*}
    Thus, if $\alpha$ is large enough, it follows that $\langle F(u_k),u_k-v \rangle$ is a negative constant for all $k$, and thus $F$ is \emph{not} Ky--Fan hemicontinuous.
\end{example}

An interesting question that remains is where the argumentation from \cite{Sadeqi2015} is incorrect. The following is a particular error which is contained in that paper.

\begin{remark}[An error in \cite{Sadeqi2015}]\label{Rem:ErrorInReference}
    Proposition~3.5 of \cite{Sadeqi2015} is wrong. Choose $K:=E:=\ell^2$ as the space of square-summable real sequences, $F:=\operatorname{Id}_E$ the identity mapping on $E$, and $\{u_k\}$ the sequence of unit vectors. Then $u_k\wto u:=0$ and
\begin{equation*}
    \langle F(u_k),u_k-u \rangle =
    \langle u_k, u_k \rangle =1
    \quad\text{for all }k.
\end{equation*}
    Hence, the sequence $\{ \langle F(u_k),u_k-u \rangle \}$ is convergent, but its limit is not equal to zero.
\end{remark}

\section*{Acknowledgements}

The author would like to thank Daniel Wachsmuth for the basic idea of the counterexample, and Ildar Sadeqi for the productive discussion leading to the creation of this paper.

\bibliographystyle{abbrv}
\bibliography{PseudoKyFan}

\end{document}